\documentclass{amsart}
 \usepackage{amsmath}

\newtheorem{theorem}{Theorem}
\newtheorem{corollary}{Corollary}
\newtheorem{lemma}{Lemma}
\newtheorem{proposition}{Proposition}

\newtheorem{definition}{Definition}

\newtheorem{example}{Example}

\newcommand{\bF}{{\mathbb F}}
\newcommand{\bZ}{{\mathbb Z}}
\newcommand{\bP}{{\mathbb P}}

\title[Ubiquity of Order Domains]{The Ubiquity of Order Domains for the Construction of Error Control Codes}

\author{John B. Little}
\address{Department of Mathematics and Computer Science\\
College of the Holy Cross\\
Worcester, MA 01610, USA}

\email{little@mathcs.holycross.edu}

\thanks{Research carried out at MSRI supported in part by 
NSF grant DMS-9810361}

\begin{document}

\begin{abstract}
Order domains are a class of commutative
rings introduced by H\o holdt,
van Lint, and Pellikaan to simplify the theory
of error control codes using ideas from algebraic geometry.
The definition is largely motivated by the structures utilized 
in the Berlekamp-Massey-Sakata (BMS) decoding algorithm, with 
Feng-Rao majority voting for unknown syndromes, applied to 
one-point geometric Goppa codes constructed from curves.  
However, order domains are much more general, and O'Sullivan 
has shown that the BMS algorithm can be used to decode codes 
constructed from order domains by a suitable generalization 
of Goppa's construction for curves.  In this article we will 
first discuss the connection between order domains and valuations 
on function fields over a finite field.  Under some mild
conditions, we will see that a general projective variety 
over a finite field has projective models which can 
be used to construct order domains and Goppa-type codes 
for which the BMS algorithm is applicable.  We will then 
give a slightly different interpretation of Geil and 
Pellikaan's extrinsic characterization of order domains via 
the theory of Gr\"obner bases, and show that their results 
are related to the existence of toric deformations of 
varieties.  To illustrate the potential usefulness of 
these observations, we present a series of new explicit 
examples of order domains associated to varieties with 
many rational points over finite fields: Hermitian hypersurfaces, 
Deligne-Lusztig varieties, Grassmannians, and flag varieties.
\end{abstract}

\maketitle
\section{Introduction}

The notion of an order domain was introduced by H\o holdt,
van Lint, and Pellikaan in \cite{14} to simplify and extend 
the theory of error control codes using ideas from algebraic 
geometry.  The definition (see Definition~\ref{oddef} below 
for the formulation we will use) is largely motivated by the 
structures utilized in the Berlekamp-Massey-Sakata (BMS) decoding 
algorithm for one-point geometric Goppa codes constructed from 
curves, the Feng-Rao bound on the minimum distance for those 
codes, and the majority voting process for unknown syndromes.  
All of these coding theoretic constructions 
are based on the properties of the ring of rational functions
with poles only at one smooth $\bF_q$-rational
point $Q$ on a curve, $R = \cup_{m=0}^\infty L(mQ)$.
These rings are the prototypical examples of order domains, 
and they furnish all the examples whose fields of fractions
have transcendence degree 1 over the field of constants.  

Geil and Pellikaan (\cite{8},\cite{6}) and 
O'Sullivan (\cite{17}) have studied the structure of order domains whose 
fields of fractions have arbitrary transcendence degree. 
Moreover, O'Sullivan (\cite{18}) has shown that the Berlekamp-Massey-Sakata 
decoding algorithm (abbreviated as the BMS algorithm in the following)
and the Feng-Rao procedure extend in a natural way to suitable 
classes of evaluation and dual evaluation codes constructed
in this much more general setting.
(See also \cite{4}, Chapter 10 for an introduction to this topic.)  

Order domains can be constructed either {\it intrinsically} or 
{\it extrinsically}, that is, by means of the algebra of 
the field $K = QF(R)$, or by means of explicit presentations
(e.g. as affine algebras $\bF_q[X_1,\ldots,X_s]/I$ with 
$I$ of a special form). 

{}From the intrinsic point of view,
the most important fact is the observation exploited
by O'Sullivan that order functions come from {\it valuations} 
on $K$.  (See \cite{25} or \cite{27}, Chapter VI for general discussions 
of valuations on fields.)  This is clear in a sense from the 
definition (see Definition~\ref{oddef} below).  The examples of order domains 
in function fields of curves also make this transparent.
Indeed, let ${\mathcal X}$ be a smooth projective curve
defined over $\bF_q$, and let $Q$ be an $\bF_q$-rational point
on ${\mathcal X}$.  Then for $R = \cup_{m=0}^\infty L(mQ)$,
$\Gamma$ is equal to the Weierstrass semigroup of ${\mathcal X}$ at $Q$ 
(the sub-semigroup of $\bZ_{\ge 0}$ consisting
of all pole orders of rational functions 
on ${\mathcal X}$ with poles only at $Q$), and $\rho(f) = -v_Q(f)$, 
where $v_Q$ is the discrete valuation at $Q$ on the function field of 
${\mathcal X}$.  

O'Sullivan extends this valuation-theoretic 
point of view to the case of function fields of transcendence
degree $\ge 2$ over $\bF_q$ (function fields of surfaces and 
higher-dimensional varieties) in \cite{17}.
He shows that every function field of transcendcence
degree 2 contains order domains
of several different types, corresponding to some of the 
possible valuation rings in these fields in a complete 
classification due originally to Zariski, and reworked
in modern language by Spivakovsky (see \cite{26} and \cite{23}).  

The valuation-theoretic interpretation of order domains
also makes connections with earlier work of Sweedler, \cite{24},
Beckman and St\"uckrad, \cite{3}, and work 
of Mosteig and Sweedler, \cite{16}, where filtrations
of rings arising from valuations are used as the foundation
for a theory of normal forms and generalized Gr\"obner bases.
\cite{17} and \cite{8} also discuss an extension
of the theory of Gr\"obner bases to ideals in order domains.

{}From the extrinsic point of view, the following result
of Geil and Pellikaan is also extremely useful, though
it applies only in the case that the value semigroup $\Gamma$
is finitely generated, hence isomorphic to a sub-semigroup of 
${\bZ}_{\ge 0}^r$ for some $r$. 
In the following statement, $M$ is an $r\times s$ matrix with entries
in $\bZ_{\ge 0}$ with linearly independent rows.  
For $\alpha\in \bZ_{\ge 0}^s$ (written as 
a column vector), the matrix product $M\alpha$ is a vector in 
$\bZ_{\ge 0}^r$.  We will call this the $M$-weight of 
the monomial $x^\alpha$.  We write $\langle M\rangle$ for the subsemigroup 
of $\bZ_{\ge 0}^r$ generated by the {\it columns} of $M$,
ordered by any convenient monomial order $\succ$ on $\bZ_{\ge 0}^r$
(for instance the $lex$ order as in Robbiano's characterization
of monomial orders by weight matrices).
We will make use of the monomial orders $>_{M,\tau}$
on $\bF_q[X_1,\ldots,X_s]$ defined as follows:
$X^\alpha >_{M,\tau} X^\beta$ if 
$M\alpha \succ M\beta$, or if $M\alpha = M\beta$ and 
$X^\alpha >_\tau X^\beta$, where $\tau$ is another
monomial order used to break ties. 

\begin{theorem}[Geil-Pellikaan, \cite{8}]
\label{GP}
Let $\Gamma = \langle M\rangle \subset \bZ_{\ge 0}^r$ 
be a semigroup. 
\begin{enumerate}
\item
Let $I \subset \bF_q[X_1,\ldots,X_s]$ be an ideal, and let $G$  
be the reduced Gr\"obner basis for $I$ with respect to a weight
order $>_{M,\tau}$ as above (abbreviated as $>$ below).  
 Suppose that 
every element of $G$ has exactly two monomials of highest 
$M$-weight in its support, and that the monomials 
in the complement of $LT_>(I)$ (the ``standard monomials''
or monomials in the ``footprint of the ideal'')
have distinct $M$-weights.  
Then $R = \bF_q[X_1,\ldots,X_s]/I$ 
is an order domain with value semigroup $\Gamma$  
and order function $\rho$ defined as follows:
Writing $f$ in $R$ as a linear combination
of the monomials in the complement of $LT_>(I)$,
$\rho(f) = \max_\succ\{M\beta: X^\beta \in supp(f)\}$.
\item  Every order domain with finitely-generated semigroup 
$\Gamma = \langle M\rangle$ has a presentation
$R \cong \bF_q[X_1,\ldots,X_s]/I$ such that the reduced
Gr\"obner basis of $I$ with respect to $>_{M,\tau}$ and 
the standard monomials are as in part 1.
\end{enumerate}
\end{theorem}

Our main goals in this article are to begin to indicate just how 
general the order domain construction is, and 
to show how the intrinsic and extrinsic 
characterizations of order domains can be used to construct 
codes from a number of interesting classes of higher-dimensional
algebraic varieties.

After some preliminaries on order functions and valuations in \S 2, 
we will begin in \S 3 by proving some general results
on the relation between order domains and valuation rings
in function fields, along the lines of \cite{16}.  
We will discuss several types of 
valuations on general function fields that are suitable for the 
construction of order domains, extending O'Sullivan's work for 
the case of the function fields of surfaces from \cite{17}.  
We will concentrate mainly on identifying the order
domains rather than on describing the corresponding valuation
rings, as is done in \cite{17}.

In rough terms, we will show that the function 
field of any projective variety ${\mathcal X}$ over a finite
field $\bF_q$ satisfying some relatively mild
conditions (for instance, the existence of a collection
of suitable subvarieties of ${\mathcal X}$ defined over $\bF_q$)
contains order domains $R$ of several different types. 
See Theorems \ref{ofval} and \ref{cdof} below for more precise, detailed 
statements.  

We have chosen to concentrate on the 
cases that lead to order domains with finitely generated value 
semigroups, since these are the ones likely to be of most
interest in coding theory.  We note, though, that
order domains with value semigroups that are not 
finitely generated (see \cite{17}, \S 5, and \cite{8}, Example 
9.6) will also exist in these function fields.

Because of the possibility of blowing up subvarieties
in varieties of dimension $\ge 2$ (see for instance
Chapter II, \S 7 of \cite{12}), the theory of valuations
on function fields of transcendence degree $\ge 2$ is 
necessarily significantly more subtle than in the case of 
transcendence degree 1.  In particular, to describe a certain 
valuation on the function field of a variety ${\mathcal X}$ of 
dimension $\ge 2$, it may be necessary to pass to a blow-up 
of ${\mathcal X}$.  \cite{17} discusses this in detail for valuations
related to monomial orders on $\bF_q[X,Y]$.  In order to 
keep the prerequisites in birational geometry to a minimum,
though, we will concentrate on the case where valuations
and order domains can be described directly from a given
projective model ${\mathcal X}$, without any blow-ups.  This
will suffice for our applications.

Following this, in \S 4, we will turn to the extrinsic 
approach and study the relation of order domains
with the toric (monomial) algebras that
appear as coordinate rings of affine toric varieties.  
The connection is that the conditions in Geil and Pellikaan's 
theorem are equivalent to saying that the order domain $R$ has 
a flat deformation to a toric algebra.  See Theorem~\ref{fldef} 
below for a more precise statement.  This point of view shows 
that there are intriguing connections between order domains and 
techniques of current interest in combinatorics, the theory of 
singularities, mirror symmetry, and other areas.

Next, in \S 5 and \S 6 we will present several
new explicit examples of order domains obtained from 
varieties such as Hermitian hypersurfaces of arbitrary
dimension, Grassmannians, and flag varieties.
These varieties have been studied with other tools in
this context by S. Hansen in \cite{11},
Rodier in \cite{19},\cite{20}  (also see the
references in those papers for earlier work).  
These particular varieties are interesting in this connection 
because they are examples of higher-dimensional varieties with 
large numbers of rational points over finite fields.  By our
results, they can be used to construct long codes for which 
the BMS algorithm, the Feng-Rao bound, and majority voting
for unknown syndromes apply.

The treatment of Hermitian hypersurfaces in \S 5 will
use the intrinsic approach to construct valuations 
of one of the types studied in \S 3.  We will 
then produce presentations of the corresponding order
domains as in Geil and Pellikaan's theorem (Theorem~\ref{GP}).  
A different construction based on the result of \S 4
will yield a second class of order domains whose
properties even more closely parallel the order domains from 
Hermitian curves.

In \S 6, on the other hand, we will 
study order domains from Grassmannians and flag 
varieties via the extrinsic approach.
Our results here depend on work of Sturmfels 
(\cite{22}) and Gonciulea and Lakshmibai (\cite{9})
establishing the existence of toric deformations of these 
varieties.  These examples can also be treated using 
the theory of Hodge algebras, or algebras with straightening
laws (ASL); see \cite{5}.  We note that these algebras
also give a sort of generalization of Gr\"obner basis
theory.  While there is a large overlap, the classes of 
Hodge algebras and order domains are distinct; neither
class contains the other.

The author would like to thank Ed Mosteig and Mike O'Sullivan
for enlightening conversations.
Work reported in this article was begun while 
the author was visiting MSRI during the Commutative Algebra program
in Spring 2003.   

\section{Preliminaries on Order Domains and Valuations}

Essentially following \cite{8}, we will use the following formulation 
of the definition of order domains.

\begin{definition}
\label{oddef}
Let $R$ be a $\bF_q$-algebra and
let $(\Gamma,+,\succ)$ be a well-ordered commutative semigroup.
An {\it order function} on $R$ is a surjective mapping 
$\rho : R \to \{-\infty\} \cup \Gamma$
satisfying:
\begin{enumerate}
\item $\rho(f) = -\infty \Leftrightarrow f = 0$,
\item $\rho(cf) = \rho(f)$ for all $f\in R$, all
$c \ne 0$ in $\bF_q$,
\item $\rho(f+g) \preceq \max_\succ \{\rho(f),\rho(g)\}$ for all 
$f,g \in R$,
\item if $\rho(f) = \rho(g) \ne -\infty$, then there exists
$c \ne 0$ in $\bF_q$ such that $\rho(f) \prec \rho(f-cg)$, and
\item $\rho(fg) = \rho(f) + \rho(g)$.
\end{enumerate}
\noindent
We will call $\Gamma$ the {\it value semigroup} of $\rho$.
\end{definition}

Axioms 1 and 5 in this definition imply that $R$ is an
integral domain.  In many cases, we will see 
that a ring $R$ with one order function has many others besides.  
For this reason an {\it order domain} is formally defined as a 
pair $(R,\rho)$ where $R$ is an $\bF_q$-algebra and $\rho$ is an 
order function on $R$.  However, we will only use one particular 
order function on $R$ at any one time.  Hence we will often omit 
it in referring to the order domain, and we will refer to $\Gamma$ 
as the value semigroup of $R$.  

Let $\alpha \in \Gamma$ be arbitrary.
The subsets $R_\alpha = \{f\in R: \rho(f) \le \alpha\}$
or $R_{< \alpha} = \{f\in R: \rho(f) < \alpha\}$
form filtrations of $R$ by $\bF_q$-vector subspaces.  
Axiom 4 implies that for each $\alpha$, $R_\alpha/R_{<\alpha}$
is a one-dimensional $\bF_q$-vector space. 
The terminology ``order function'' is supposed to 
suggest the existence of $\bF_q$-bases of $R$ 
whose elements have distinct $\rho$-values, and are
hence ordered by $\rho$.  This is a consequence of 
the one-dimensional quotients axiom 4.

\cite{8} and \cite{18} contain a number of
examples of order domains; we will provide additional
examples in \S 5 and \S 6.

At this point a comment concerning the relation of this
definition to the one used in \cite{17} and \cite{18}
is probably in order.  In those papers an order function
is defined as a mapping $o : R \to {\mathbb N} \cup \{-1\}$
satisfying the properties that for all $a$,
the set $L_a = \{f \in R : o(f) \le a\}$ is an $\bF_q$-vector
space of dimension  $a+1$, and $o(f) < o(g)$ implies
$o(fz) < o(gz)$ for all $f,g,z\in R$.  It can be seen that this
formulation satisfies all the axioms in Definition~\ref{oddef}, 
but it is {\it less general} than our definition.  It 
excludes, for instance, $R$ such as $\bF_q[X,Y]$ with order 
function induced by a lexicographic monomial order (that is,
the order function $\rho(X^nY^m) = (m,n) \in \Gamma = \bZ_{\ge 0}^2$,
ordered lexicographically). Note for instance that 
the lexicographic order does not satisfy Proposition 1.2
of \cite{17}.  Nevertheless, lexicographic and similar 
orders on polynomial rings do furnish examples of
order domains as in our definition.  In particular the 
well-ordering property does hold, even though there is 
no power $Y^n$ satisfying $\rho(Y^n) > \rho(X)$. 

We will follow the notation and terminology of
\cite{25} for Krull valuations on function fields.
Let $K$ be a field.  A valuation $v$ of $K$ is a mapping 
from $K$ to $\Lambda \cup \{+\infty\}$, where
$\Lambda$ is a totally ordered abelian group satisfying
\begin{enumerate}
\item $v(f) = +\infty$ if and only if $f = 0$,
\item $v(fg) = v(f) + v(g)$ for all $f,g\in K$,
\item $v(f+g) \ge \min \{v(f),v(g)\}$ for all $f,g \in K$.
\end{enumerate}
Given any valuation on $K$, the corresponding {\it valuation ring} is
$S_v = \{f \in K : v(f) \ge 0\}$,
a local ring with maximal ideal 
$M_v = \{f \in K : v(f) > 0\}$.
The {\it residue field} of $v$ is the quotient $k_v = S_v/M_v$.

We will always consider function fields $K$ with a 
constant subfield $k$ equal to a finite field $\bF_q$.
All valuations will be trivial on $k$ (i.e. $v(c) = 0$,
if $c \in k$).
The {\it dimension} of $v$ is the transcendence degree of 
$k_v$ over $k$.  We will be concerned only with 
valuations of dimension 0, so the residue field will 
be at most an algebraic extension of the constant field.

Let $\Lambda$ be the value group of a valuation.
A subset $\Sigma \subseteq \Lambda$ is said to be 
a {\it segment} if whenever $\beta \in \Lambda$ is between 
$-\sigma$ and $\sigma$ (in the order) for some 
$\sigma\in \Sigma$, then $\beta\in \Sigma$.
An {\it isolated subgroup} of $\Lambda$ is a proper 
subgroup that is also a segment.  The {\it rank} of 
a valuation $v$ is the number of isolated subgroups 
of the value group (or $\infty$ if that number is not 
finite).  

The {\it rational rank} of $v$ is the dimension of 
the ${\mathbb Q}$-vector space $\Lambda \otimes_{\bZ} {\mathbb Q}$.

The valuation $v$ is said to be {\it discrete} if its
value group is a discrete group of finite rank, that is,
isomorphic to a subgroup of ${\bZ}^n$ for some $n$, 
ordered by the lexicographic order.
 
Let $K$ be the function field of a
variety ${\mathcal X}$.
We say a valuation $v$ is {\it centered} at a 
(closed) point $Q\in {\mathcal X}$ if 
we have the containment of local rings
${\mathcal O}_{{\mathcal X},Q} \subseteq S_v$,
and the maximal ideals satisfy 
$M_v \cap {\mathcal O}_{{\mathcal X},Q} = M_{{\mathcal X},Q}$.  

\section{Constructing Order Domains from Valuations}

As shown in \cite{17} and \cite{8}, 
every order function $\rho$ on $R$ determines a 
valuation $v$ on $K = QF(R)$, defined by:
\begin{equation}
\label{val}
v(f/g) = \rho(g) - \rho(f).
\end{equation}
Note the signs; it follows that $\rho = -v|_R$.  
The {\it value group} of this $v$ is the group of differences 
$\Lambda = \Gamma - \Gamma$, and the ordering is induced
by the ordering on $\Gamma$.

{}From \eqref{val}, it can be seen that the order domain $R$ and 
the valuation ring $S_v$ for the valuation corresponding to $\rho$ 
are in a special relative position in $K = QF(R) = QF(S_v)$.
Namely, we have the following statement.

\begin{proposition}
\label{comp}
Let $R$ be an order domain,
and $S_v$ be the corresponding valuation ring of $K = R$ as above.
Then $S_v = (R\cap S_v) + M_v$ and $R \cap S_v = \bF_q$.
\end{proposition}

\begin{proof}
Both claims follow from the definitions.  We have
$S_v = \{f/g : f,g \in R, v(f) - v(g)\ge 0\}$ 
and $M_v = \{f/g : f,g \in R, v(f) - v(g) > 0\}$, while
$v(h) = -\rho(h) \le 0$ for all $h \in R$.  So 
$R \cap S_v = \{h\in R: \rho(h) = 0\}$.
By axioms 2 and 5 in Definition~\ref{oddef}, if $\rho(h) = 0$, then
$h$ must be constant (an element of $\bF_q$), and 
conversely.  So $R\cap S_v = \bF_q$, and $S_v = \bF_q + M_v$.
\end{proof}

In other words, the valuation $v$ has dimension $\dim(v) = 0$,
according to the terminology from \S 2, since
$S_v/M_v = \bF_q$.  Moreover,
the valuation ring $S_v$ is an {\it $\bF_q$-complementary 
valuation ring} to $R$ in the terminology of \cite{16}.
(Note that the definition of a valuation considered in that 
article is essentially the extension of our order function
$\rho$ to $K$:  $v(f/g) = \rho(f) - \rho(g)$, rather
than the negative as in \eqref{val}.)
In \cite{16}, the connection between valuations and filtrations
is considered in detail, and in Theorem 4.4, Lemmas 4.7 and 4.8, and 
Proposition 4.9 it is shown that there is a
one-to-one correspondence between $\bF_q$-complementary
valuation rings to a subring $R$ in $K$ and 
{\it regular, normalized filtrations} on $R$. The properties
of these regular, normalized filtrations, in particular
the one-dimensionality over $\bF_q$ of the graded quotients,
are equivalent to the axioms defining an order domain as
in Definition~\ref{oddef}. 

We now want to turn the tables and show how to produce
examples of order domains starting from valuations.
The following theorem gives a general class of 
valuations on function fields $K = K({\mathcal X})$ for which 
there are corresponding order domains $R \subset K$.
We restrict to $R$ of this special form because this
seems to be the most important case for applications
in coding theory.

\begin{theorem}
\label{ofval}
Let $R$ be an affine domain over $\bF_q$, that is 
$$R \cong \bF_q[X_1,\ldots,X_s]/I$$
where $I$ is a prime ideal.  Let ${\mathcal X}$ be the 
projective closure of $V(I)$ in ${\mathbb P}^s$, and let
$H_0$ (with reduced scheme structure) 
be the intersection of ${\mathcal X}$ with the
hyperplane at infinity.  Assume $H_0$ is an irreducible 
divisor on ${\mathcal X}$.  Let $v$ be any valuation 
on the function field $K({\mathcal X})$ such that 
\begin{enumerate}
\item the rational rank of $v$ is $d = \dim {\mathcal X}$, and 
\item $v$ is centered at a smooth point $Q\in H_0\subset {\mathcal X}$.
\item $v(f) \le 0$ in $\Lambda$ for all $f \in R$.  
\end{enumerate}
Then $\rho = -v|_R$ is an order function on $R$.
\end{theorem}

(Note that we can view $R$ as a subring of $K$ consisting 
of functions with poles on $H_0$.)

\begin{proof}
Hypothesis 1 implies that the Abhyankar
inequality for $v$:
$$rat. rank(v) + tr. deg._{\bF_q}(k_v) \le \dim({\mathcal X})$$
(see \cite{25}, Th\'eor\`eme 9.2) is 
an equality.  Hence, by part b of that theorem, 
the valuation group of $v$ is isomorphic (as a group) to $\bZ^d$. 
The rank of $v$ may be any integer $r$ with $1 \le r \le d$, 
though, so the ordering may be any one of a number
of different possibilities.
For example, we may have $\Lambda$ discrete (the case $r = d$),
$\Lambda$ a subgroup of ${\mathbb R}$ generated by $d$ ${\mathbb Q}$-linearly
independent real numbers (the case $r = 1$), or an intermediate
case.

Axioms 1, 2, 3, and 5 in the definition of 
an order function follow immediately from the definition
of a valuation, and show that $\Gamma = \rho(R)$ 
is a semigroup contained in the value group $\Lambda$ of $v$.

To show that the one-dimensional quotients 
axiom 4 holds, we will use Theorem 4.4, vii
from \cite{16}.  We must show that $S_v \cap R = \bF_q$
(that is, that $S_v$ and $R$ are in $\bF_q$-complementary
position in $K$).  This follows from the irreducibility 
of $H_0$ and hypothesis 2.
The nonconstant $f \in R$ are not contained in $S_v$
by hypothesis 3.  

What remains to be proved is that $\Gamma = \rho(R) = -v(R)$
is well-ordered.  We will use the following criterion.

\begin{lemma} 
Let $(\Gamma,+)$ be any finitely
generated inverse-free semigroup.  If $\prec$ is any
total order on $\Gamma$ compatible with the addition
operation (in the sense that $\alpha \prec \beta$ implies
$\alpha + \gamma \prec \beta + \gamma$ for all $\alpha, \beta,
\gamma\in \Gamma$), then $\Gamma$ is well-ordered under 
$\prec$.
\end{lemma}

A proof of the Lemma is sketched on page 371 of \cite{8}.

Our $\Gamma$ is clearly inverse-free, since $\Gamma$ is 
contained in the set of elements of the value group 
$\Lambda$ that are $\ge 0$.  The order on $\Gamma$ is
induced from that on $\Lambda$, so is compatible with 
addition.  So we are reduced to showing that
$\Gamma$ is finitely generated.  This follows, for instance,
from the Noether Normalization theorem (see \cite{27}, Chapter VII,
\S 7, Theorem 35, or \cite{10}, Theorem 3.4.1, and Exercises 3.4.1
and 3.4.2).  There exists a transcendence
basis $\{z_1,\ldots,z_d\}$ in $R$ such that $R$ is a finite,
integral extension of the polynomial ring 
$\bF_q[z_1,\ldots,z_d]$.  $\Gamma$ is generated by 
the values $\rho(z_i)$ and the $\rho$-values for
the elements of a basis of $R$ over $\bF_q[z_1,\ldots,z_d]$.
Indeed, we obtain from the Gr\"obner basis algorithm
for Noether Normalization described in \cite{10} a
monomial $\bF_q$-basis for $R$ consisting of products
of arbitrary monomials in the $z_1,\ldots,z_d$ with
a finite list of monomials in the remaining variables.
The values of $\rho$ on these basis monomials are distinct.
Otherwise, we would have an algebraic dependence because
of the one-dimensional quotients property.
\end{proof}

The $\rho$ given by this theorem 
are all {\it monomial} order functions as defined
in \cite{17}.  To construct one class of valuations as in Theorem~\ref{ofval}
in a simple fashion, starting from the geometry 
of the variety ${\mathcal X}$, we can use the well-known
{\it composite divisorial valuations} described,
for instance, in \cite{25}, Example 9 and the following 
remark (see also \cite{2} for relations 
between these valuations and monomial orderings 
in the theory of Gr\"obner bases).  

Let  ${\mathcal X}$ be a projective variety of dimension $d$ defined
over a finite field $\bF_q$, and let 
\begin{equation}
\label{flag}
{\mathcal F}: {\mathcal X} = V_0 \supset V_1 \supset V_2 \supset 
\cdots \supset V_d
\end{equation}
be a flag of subvarieties of ${\mathcal X}$ satisfying the following 
conditions:
\begin{enumerate}
\item Each $V_i$ is irreducible and defined over $\bF_q$.
\item The dimension of $V_i$ is $r - i$ for each $i$ (so $i$ is
the codimension in ${\mathcal X}$).
\item For each $0 \le i\le d-1$, $V_i$ is smooth at the generic point 
of $V_{i+1}$.
\end{enumerate}

Since each $V_{i+1}$ is an irreducible divisor in $V_i$, 
each rational function $g$ on $V_{i+1}$ has a well-defined {\it (vanishing
or pole) order} along $V_{i+1}$, denoted $ord_{V_{i+1}}(g)$.  
By definition, $ord_{V_{i+1}}(g)$ is positive
if $g$ vanishes along $V_{i+1}$, and negative if $g$ has a pole
along $V_{i+1}$.  
We note that it also follows from these hypotheses that $V_d$ is a smooth
$\bF_q$-rational point on the irreducible curve $V_{d-1}$.

Any such flag ${\mathcal F}$ defines a valuation $v_{\mathcal F}$ on 
the function field $K = K({\mathcal X})$ as follows.  For each $i$, 
fix some function $g_i$ on the subvariety $V_{i-1}$ 
with a zero of order 1 along $V_i$.  
Given any $f \in K$, we define a sequence of integers 
(the notation $F|_{V_i}$ means 
the function $F$, restricted to the variety $V_i$)
\begin{align*}
v_1 &= ord_{V_1}(f)\\
v_2 &= ord_{V_2}\left((f/g_1^{v_1})|_{V_1}\right)\\
    &\ \,\vdots& \\
v_d &= 
ord_{V_d}\left((f/(g_1^{v_1}\cdots g_{d-1}^{v_{d-1}}))|_{V_{d-1}}\right),
\end{align*}
and let
\begin{equation}
\label{flagval} 
v_{\mathcal F}(f) = (v_1,\ldots,v_d) \in \bZ^d.
\end{equation}

Then $v_{\mathcal F}$ is a discrete valuation of $K$ with 
rational rank $d$, rank $d$, and value group 
$\bZ^d$, ordered lexicographically.  The values
$v_{\mathcal F}(f)$ depend on the choice of the
auxiliary functions $g_i$.  However all choices 
of $g_i$ will lead to equivalent orderings.  Indeed, note that the
auxiliary functions are unnecessary for {\it comparing
valuations} of two functions $f$ and $f'$; the comparison
can be made using only the orders $ord_{V_1}(f/f')$, 
$ord_{V_2}((f/f')|_{V_1})$ and so on.

For example, it is easy to see that the lexicographic
order on $R = \bF_q[X_1,\ldots,X_s]$ (the affine
coordinate ring of a standard affine subset of 
${\mathcal X} = {\mathbb P}^s$)
with $X_1 > X_2 > \cdots > X_s$ is obtained
from the composite divisorial valuation
$v_{\mathcal F}$ with 
$${\mathcal F}: {\mathbb P}^s \supset V(X_1) \supset V(X_1,X_2) 
\supset \cdots \supset V(X_1,\ldots,X_s).$$
The center is the origin in the affine plane, and this
shows that lexicographic order functions on the polynomial ring $R$
do not come from the construction of Theorem~\ref{ofval}.
Similarly, it is not difficult to see that 
graded reverse lexicographic monomial orders on
this $R$ come from composite divisorial valuations
constructed from subvarieties on a blow-up of ${\mathbb P}^s$
(see \cite{2}, \S 1, for example).
In the following, if we wanted to consider order functions
like the graded reverse lexicographic order, in effect,
${\mathcal X}$ would be the blow-up of ${\mathbb P}^s$.

\begin{theorem}
\label{cdof}
Let ${\mathcal X}$ be any projective
variety over $\bF_q$ which has a flag ${\mathcal F}$ of subvarieties
defined over $\bF_q$ satisfying the hypotheses above, and
such that $V_1$ is ample on ${\mathcal X}$ (that is, 
the complete linear system $|\ell V_1|$ defines
a projective embedding of ${\mathcal X}$ for some $\ell \ge 1$). 
Let $v_{\mathcal F}$ be the corresponding valuation of
the function field $K = K({\mathcal X})$ defined in \eqref{flagval}.
Let $R$ be the subring $R = \cup_{m=0}^\infty L(m V_1)$
(the subring of $K$ consisting of functions
with poles only along the subvariety $V_1$).  Let
$$\rho(f) = -v_{\mathcal F}|_R(f)$$
if $f \ne 0$, and $\rho(f) = -\infty$ if $f = 0$.
Then $\rho$ is an order function on $R$. 
\end{theorem}

\begin{proof}
This follows from Theorem~\ref{ofval} on 
reembedding ${\mathcal X}$ so that the divisor $\ell V_1$ becomes
the hyperplane section $H_0$.  The center
of the valuation is the point $Q = V_d$.
The $v_{\mathcal F}$ valuations have rational rank
$d = \dim {\mathcal X}$ by construction.

The well-ordering property of the image of $\rho$
also follows by a direct argument in this case.  
Suppose we had an infinite strictly descending chain:
\begin{equation}
\label{desc}
\rho(f_1) > \rho(f_2) > \rho(f_3) > \cdots,
\end{equation}
where $f_i \in R$ for all $i$ and $>$ denotes the lexicographic 
order in $\bZ^r$.  The $f_i$ must be non-constant,
so by the definition of $R$ the first components 
of $\rho(f_i)$ (that is, the integers $-ord_{V_1}(f_i)$) 
are strictly positive.  This follows since there are
no nonconstant functions in $R$ with a pole of order $\le 0$
along $V_1$ (recall that functions in $R$ can have poles only 
along $V_1$).  Hence the first components stabilize after a 
finite number of steps in the chain \eqref{desc}: There exists
$i_0 \ge 1$ such that $ord_{V_1}(f_j) = ord_{V_1}(f_{i_0}) = n$
for all $j \ge i_0$ and some $n \ge 1$.  The set of rational
functions on ${\mathcal X}$ with a pole of order at most $n$ 
along $V_1$ and no other poles (together with the zero function)
forms a finite-dimensional vector space over the field
of constants (this follows, for example from \cite{12},
Chapter II, Theorem 5.19).  
Hence the orders of poles and zeroes of the 
$f_j/(g_1^{n}\cdots g_{k-1}^{v_{k-1}})|_{V_k}$, $j \ge i_0$
along the $V_k$, $k \ge 2$ are bounded.  As a result, the 
chain \eqref{desc} must eventually stabilize.  (It would 
also be possible to find a sub-semigroup of $\bZ_{\ge 0}^d$ 
isomorphic to the image of $\rho$. An example of this is 
given in \S 4 below.) 
\end{proof}

Note that when $r \ge 2$, the choice of $V_1$ determines 
the ring $R$, but the choice of rest of the flag ${\mathcal F}$
still possibly yields many different order functions on $R$.  Moreover, 
the rest of the flag is necessary because the pole order
$ord_{V_1}(f)$ alone gives a filtration of $K$ too coarse
to satisfy the one-dimensional quotients axiom 4 in the 
definition of an order function.  As noted before, $V_d$ 
is a smooth point of the irreducible curve $V_{d-1}$.  
Hence $ord_{V_d}$ is a discrete, rank 1 valuation on the 
function field of $V_{d-1}$ and the one-dimensional quotients 
property in axiom 4 also follows directly for these valuations 
from this observation.  

Even if a given variety ${\mathcal X}$ defined
over $\bF_q$ does not have any suitable flags of subvarieties
defined over the field $\bF_q$, they always exist over the 
algebraic closure $\overline{\bF_q}$, hence over some finite 
extension of $\bF_q$.  

Another way to frame what Theorem~\ref{cdof} says is the following.
If we think of the divisor $V_1$ in the flag used in the proof 
as the intersection of the image of ${\mathcal X}$ under this new 
embedding with the hyperplane at infinity, then the coordinate 
ring of the affine variety ${\mathcal X}\ \backslash V_1$ will 
have order functions. 

For a very simple example, consider the ring 
$R = \bF_q[X,Y]/\langle XY-1\rangle$ cited in \cite{14} 
as a ring with no order functions as in Definition~\ref{oddef}.  
In geometric terms, this is the coordinate ring of an 
affine subset of the projective conic $V(XY - Z^2)$, which 
is birationally isomorphic to ${\mathcal X} = {\mathbb P}^1$ 
(the function field is isomorphic to the field of rational 
functions in one variable).  There are many ways to embed 
${\mathbb P}^1$ to yield examples of order domains. In fact
in this case, a simple linear change of coordinates to make 
the hyperplane at infinity tangent to the curve will produce
a conic curve that does yield an order domain:  
$R' = \bF_q[X,Y]/\langle X - Y^2\rangle$.

\begin{corollary}
Any projective 
variety ${\mathcal X}$ which has a suitable
flag of subvarieties defined over $\bF_q$ can be used to construct 
evaluation and dual evaluation codes for which the BMS algorithm 
and the Feng-Rao procedure are applicable.
\end{corollary}

Though we have considerable
freedom in the choice of the flag ${\mathcal F}$, the
composite divisorial valuations 
$v_{\mathcal F}$ in no way exhaust the valuations described in 
Theorem~\ref{ofval}.

\section{Order Domains from Toric Deformations}

This section is devoted to the observation that the extrinsic 
characterization of order domains from Theorem~\ref{GP} can be 
reinterpreted in another way giving a way to generate 
additional examples of order domains.  

Given a semigroup $\Gamma \subseteq {\bZ}_{\ge 0}^r$,
let $R_\Gamma$ be the subring of $\bF_q[T_1,\ldots,T_r]$ generated by
the monomials $T^\gamma$ for $\gamma \in \Gamma$. 
By \cite{8}, Proposition 4.8, $R_\Gamma$ is
an order domain with value semgroup $\Gamma$. 
Moreover, the graded algebra $Gr(R)$ of any order
domain $R$ with value semigroup $\Gamma$ is 
isomorphic to $R_\Gamma$ by \cite{8}, Proposition 6.5.

A first connection
with the extrinsic characterization in Theorem~\ref{GP} is given
by \cite{8}, Proposition 10.6.
Let $M$ be the $r\times s$ matrix whose columns are
$\gamma_1, \ldots, \gamma_s$ in some generating set for
$\Gamma$.  We define $I_\Gamma$ to be the 
binomial ideal generated by the $X^\alpha - X^\beta$ such 
that $M\alpha = M\beta$.  In the literature, 
$I_\Gamma$ is also known as the {\it toric ideal} 
corresponding to $\Gamma$.  The variety $V(I_\Gamma)$ 
is toric (but is not necessarily a normal variety, as
is sometimes required in the study of toric varieties).

In this section we will prove the following complement
of these results.

\begin{theorem} 
\label{fldef}
Let $R$ be an order domain with a given 
finitely-generated value semigroup
$\Gamma \subset \bZ_{\ge 0}^r$.  Let 
$$R_{\Gamma} = \bF_q[\Gamma] \cong \bF_q[X_1,\ldots,X_s]/I_\Gamma$$
as above.  Then $R$ has a flat deformation to $R_\Gamma$. 
Conversely, every flat deformation of $R_{\Gamma}$ of the form 
given in Theorem~\ref{GP}
is an order domain with value semigroup $\Gamma$.
\end{theorem}

\begin{proof}
Let $R$ be an order domain with value semigroup $\Gamma$.
By part 2 of Geil and Pellikaan's theorem, 
we have a presentation 
$$R \cong \bF_q[X_1,\ldots,X_s]/I,$$ where 
$I$ has a Gr\"obner basis of the form described in part 1 
of Theorem~\ref{GP}.  The deformation can be seen explicitly as follows.
Let $\omega$ be a general linear combination
of the rows of the matrix $M$ in Theorem~\ref{GP}.  Then we get
a one-parameter family of varieties over ${\mathbb A}^1_t$
by mapping $X_i \to t^{-\omega_i} X_i$ in the 
Gr\"obner basis elements for $I$.  Clearing denominators,
the terms of non-maximal $M$-weight in the generators of $I$
vanish on letting $t\to 0$, and the limiting ideal is the binomial 
ideal $I_\Gamma$.  Flatness follows from the requirement in 
Theorem~\ref{GP} that the 
specified generators of $I$ form a Gr\"obner basis. (See, for 
instance, \cite{10}, section 7.5.  In our case, the deformation
is to the binomial ideal $I_\Gamma$ rather than a monomial ideal,
but the idea is the same.)
The converse is just a restatement of Geil and Pellikaan's theorem.
\end{proof}

The corresponding affine varieties $V(I)$ defined 
by the ideals $I$ in Theorem~\ref{GP}
are flat deformations of the toric variety $V(I_\Gamma)$.

In the next sections, we will consider several examples of
the way order domains can be identified from explicit
varieties of interest in coding theory.
We will show the existence of order functions both by 
using Theorem~\ref{cdof} and by using Theorem~\ref{fldef}.

\section{Order Domains and Codes from Hermitian Hypersurfaces}

In this section, we will begin by considering order domains associated
to Hermitian hypersurfaces in $\bP^{r+1}$ over a field
$\bF_{q^2}$, for any $r\ge 1$.  These varieties have a number
of properties that make them interesting for the construction
of codes, such as large numbers of $\bF_{q^2}$-rational
points, large automorphism groups, and so forth.
The case of Hermitian curves in
the plane has, of course, been extensively studied by many 
coding theorists over the past 15 years.  Hermitian hypersurfaces
have been considered by S. Hansen (\cite{11}) and 
Rodier (\cite{19}). 

We will take the following as the standard projective embedding
of the Hermitian hypersurface:
\begin{equation}
\label{herm}
{\mathcal H}_r = V(X_0^{q+1} + X_1^{q+1} + \cdots + 
X_r^{q+1} + X_{r+1}^{q+1}) \subset {\mathbb P}^{r+1}.
\end{equation}

\begin{example}
We will begin by constructing an explicit order domain
from the Hermitian surface ${\mathcal H}_2$ in ${\mathbb P}^3$.  
If we use the standard form \eqref{herm} for the defining equation,
then the hyperplane section $V_1 = V(X_0) \cap {\mathcal H}_2$ is a
smooth Hermitian curve in the plane $V(X_0)$.  We can then
select any $\bF_{q^2}$-rational point on $V_1$ to be the point
$V_2$ in a flag ${\mathcal F} : {\mathcal H}_2 \supset V_1 \supset V_2$
as in \S 2.  For example, consider $V_2 = (0:1:\delta:0)$,
where $\delta^{q+1} = -1$ in $\bF_{q^2}$. The corresponding ring 
$R = \cup_{m=0}^\infty L(mV_1)$
thus has the structure of an order domain, where $\rho$ is 
constructed from a composite divisorial valuation as in Theorem~\ref{cdof}.  

To work explicitly with 
$R$, note that the rational functions that are 
contained in $R$ (in which the denominators 
can only contain powers of $X_0$) can be identified with polynomials 
in $X_1,X_2,X_3$, after dehomogenizing with respect to $X_0$.
We have the following results by easy computations:
\begin{align*}
           \rho(X_1) &= (1,0)\\
           \rho(X_2) &= (1,0)\\
           \rho(X_3) &= (1,-1)
\end{align*}
Moreover, by axiom 4 in Definition~\ref{oddef}, we expect some 
linear combination
of $X_1$ and $X_2$ to have smaller $\rho$-value than $\rho(X_1) = \rho(X_2)$.
This linear combination comes from the equation of the tangent line
to the curve $V_1$ at the point $V_2$:
$$\rho(\delta X_1 - X_2) = \rho(X_1 + \delta^q X_2) = (1,-(q+1)).$$
We will write $U = \delta X_1 - X_2$ in the following.

Since there are no linear polynomials in $X_1,X_2,X_3$ on $V_1$ 
vanishing to order higher than $q+1$ at $V_2$, we can
also obtain an order-preserving linear mapping from the value group
of $\rho$ to a sub-semigroup of $\bZ_{\ge 0}^2$ (with the lexicographic
order) by mapping $(a,b) \mapsto (a,(q+1)a+b)$.  Composing with $\rho$ gives
a new order function $\tilde{\rho} : R \to \Gamma \subset \bZ_{\ge 0}^2$
as in Geil and Pellikaan's extrinsic characterization of 
order domains.  We have
\begin{align}
\label{vals}
           \tilde{\rho}(X_1) &= (1,q+1)\notag\\
           \tilde{\rho}(X_2) &= (1,q+1)\notag\\
           \tilde{\rho}(X_3) &= (1,q)\\
           \tilde{\rho}(U) &= (1,0)\notag
\end{align}

To put $R$ into Geil and Pellikaan's form, we will make a 
linear change of coordinates, substituting $X_2 = \delta X_1 - U$ and
writing the equation in terms of $X_1,U,X_3$.   (From Theorem~\ref{GP}, 
recall that the columns of the weight matrix $M$ should generate the 
value semigroup of $R$, so from \eqref{vals} we see that $U$ should 
be used to replace $X_2$.)  The result is the equation:
\begin{equation}
\label{newherm}
F(X_1,U,X_3) = X_3^{q+1} + \delta^q X_1^q U + \delta X_1 U^q - U^{q+1} - 1 = 0.
\end{equation}
We use the monomial order $>_{M,lex}$, where
$$M = \begin{pmatrix}
               1 & 1 & 1\\
             q+1 & 0 & q\\
      \end{pmatrix}$$
comes from the $\tilde{\rho}$-values above.
Note that the columns of $M$ correspond to $X_1,U,X_3$ in that
order, but the variables are ordered $X_1 > X_3 > U$
according to $>_{M,lex}$.  We see that
$F$ has exactly two terms of highest $M$-weight (with the monomials
$X_3^{q+1}$ and $X_1^q U$).  Moreover, the leading term is $X_3^{q+1}$,
and the monomials in the complement of the initial ideal are
\begin{equation}
\label{deltaset}
\Delta = \{X_1^a U^b X_3^c :  a,b \ge 0; 0 \le c \le q\}.
\end{equation}
If two monomials   $X_1^a U^b X_3^c$ and 
$X_1^{a'} U^{b'} X_3^{c'}$ in $\Delta$ have the same $M$-weight, 
then $a+b+c = a'+b'+c'$ and $(q+1)a + qc = (q+1)a'+ qc'$.
The second equation says $c - c'$ is divisible by $q+1$.  But this
is only possible if $c = c'$, and $a = a'$ and $b = b'$ follow.
Thus, we have verified the hypotheses of Geil and Pellikaan's
extrinsic characterization in Theorem~\ref{GP}, 
so we have proved that $\tilde{\rho}$ is an order function.

The Hermitian surface ${\mathcal H}_2$ has $(q^2+1)(q^3+1)$ 
$\bF_{q^2}$-rational
points, of which $q^3 + 1$ lie on $V(X_0)$.  Hence evaluation and 
dual evaluation codes of length $n = q^2(q^3+1)$ can be constructed
by this method from $R$, using the ordering induced by $\tilde{\rho}$
on the monomials in $\Delta$ from \eqref{deltaset}.   We will discuss 
the analogous codes constructed from all ${\mathcal H}_r$ in more
detail below.  The BMS algorithm applies for decoding the dual 
of the evaluation code with basis formed by evaluation 
of the first $\ell$ monomials in $\Delta$ (in the $>_{M,lex}$ 
order), for any $\ell \ge 1$.
$\diamondsuit$
\end{example}

The construction in this example can be extended by a natural
inductive procedure to define order domains from all the 
Hermitian hypersurfaces.

\begin{proposition}
Let ${\mathcal H}_r$ be the Hermitian
hypersurface in ${\mathbb P}^{r+1}$.  We can construct an order domain
from ${\mathcal H}_r$ with value semigroup in $\bZ^r$ generated by 
\begin{align*}
           &(1,0,0,0\cdots,0)\\
           &(1,-1,0,0,\cdots,0)\\
           &(1,0,-1,0,\cdots,0)\\
           &\qquad \vdots\\
           &(1,0,0,0,\cdots,-1)\\
           &(1,0,0,0,\cdots,-(q+1)),
\end{align*}
and presentation of the form
$$R\cong \bF_{q^3}[X_1,\ldots,X_{r-1},X_{r+1},U]/\langle F\rangle$$
where
$$F = X_{r+1}^{q+1} + \delta^q X_{r-1}^q U + \delta X_{r-1} U^q + 
\sum_{j=1}^{r-2} X_{j}^{q+1} - U^{q+1} - 1.$$
\end{proposition}

\begin{proof}
Let $\delta^{q+1} = -1$ in $\bF_{q^2}$.
We construct a flag ${\mathcal F}$ of subvarieties of 
${\mathcal H}_r$ as follows.  $V_0 = {\mathcal H}_r$, 
$V_i = V(X_0,\ldots,X_{i-1})\cap {\mathcal H}_r$ for $i = 1,\ldots, r-1$,
and 
$$V_r = V(X_0,\ldots,X_{r-2},\delta X_{r-1} - X_r)\cap {\mathcal H}_r.$$
It is easy to see that the hypotheses for Theorem~\ref{cdof} are 
satisfied, because
for $0\le i \le {r-1}$, $V_i$ is a smooth $(r-i)$-dimensional
Hermitian variety in the linear space $V(X_0,\ldots,X_{i-1})$, 
and $V_r = \{(0,0,\ldots,1,\delta,0)$ is a smooth point of $V_{r-1}$.
We choose the auxiliary functions $g_i = X_{i-1}/X_{r-1}$, and 
define $\rho = -v_{\mathcal F}|_R$ as in \S 3.  (We view
$V_1$ as the intersection of ${\mathcal H}_r$ with 
the hyperplane at infinity and dehomogenize
by setting $X_0 = 1$.)  Then
\begin{align*}
           \rho(X_1) &= (1,-1,0,0,\cdots,0)\\
           \rho(X_2) &= (1,0,-1,0,\cdots,0)\\
                     &\ \,\vdots \\
           \rho(X_{r-2}) &= (1,0,0,0,\cdots,-1,0)\\
           \rho(X_{r-1}) &= \rho(X_r) = (1,0,0,0\cdots,0)\\
           \rho(X_{r+1}) &= (1,0,0,0,\cdots,-1)\\
           \rho((\delta X_{r-1} - X_r)) &= (1,0,0,0,\cdots,-(q+1))
\end{align*}
(The last comes from the linear combination $U =\delta X_{r-1} - X_r$ 
defining the tangent line to $V_{r-1}$ at the point $V_r$). 
The extrinsic form of the corresponding order domain $R$
is directly parallel to \eqref{newherm}. The first two monomials
have the largest, and equal weights. 
\end{proof}

We now present some more detailed information concerning some codes
constructed from these Hermitian hypersurfaces.  So that we will 
be in the setup for O'Sullivan's generalized BMS algorithm, we
will consider codes obtained by evaluation of polynomials in 
$X_1,\ldots,X_{r-1}, U, X_{r+1}$ at the affine $\bF_{q^2}$-rational points
on ${\mathcal H}_r$.  Let $C_a$, $a\ge 1$, denote the code defined by the 
monomials of total degree $\le a$.  (Note that since the first
component of $\rho(X_j)$ is $1$ for all $j$, the orders $>_{M,\tau}$ 
are all graded orders, so if $a \le q$, the $C_a$ code is spanned by the 
evaluations of the first $\binom{a + r + 1}{a}$ monomials in this order 
and has the form needed for BMS.)

\begin{theorem}
\label{hermcodes}
The $C_1$ code over $\bF_{q^2}$ defined in this way
from ${\mathcal H}_r \subset \bP^{r+1}$ has the following parameters:
\begin{align*}
         n &= q^{2r +1} - (-1)^{r+1} q^r\\
         k &= r+2\\
         d &= \begin{cases} q^{2r+1} - q^{2r-1} & \text{\rm if $r$ is even} \\
               q^{2r+1} - q^{2r-1} - q^r - q^{r-1} & \text{\rm if $r$ is odd}
               \end{cases}
\end{align*}
\end{theorem}

\begin{proof}
We use the results from \cite{1}.  The number of $\bF_{q^2}$-rational
points of the Hermitian hypersurface ${\mathcal H}_r$ is
\begin{equation}
\label{nopts}
|{\mathcal H}_r(\bF_{q^2})| = 
\frac{(q^{r+2} + (-1)^{r+1})(q^{r+1} - (-1)^{r+1})}{q^2 - 1}
\end{equation}
The intersection of ${\mathcal H}_r$ with the hyperplane $X_0 = 0$ at 
infinity is isomorphic to the Hermitian variety ${\mathcal H}_{r-1}$.
Hence by \eqref{nopts}, the number of affine $\bF_{q^2}$-rational points
is 
$$|{\mathcal H}_r(\bF_{q^2})| - |{\mathcal H}_{r-1}(\bF_{q^2})|
= q^{2r + 1} - (-1)^{r+1} q^r.$$
This yields the blocklength $n$ of the $C_a$ codes for all $a\ge 1$.  
The dimension of $C_1$ is $k = r + 2$ since the codewords obtained
by evaluation of $1,X_1,\ldots,X_{r-1},X_{r+1},U$ are linearly 
independent.

To determine the minimum distance $d$, we must determine the 
largest possible number of {\it affine} $\bF_{q^2}$-rational
points in an intersection ${\mathcal H}_r \cap L$ where $L$
is the hyperplane in $\bP^{r+1}$ defined by a linear 
form.  By \cite{1} again, there are
exactly two cases for $L$ defined over $\bF_{q^2}$.
Either $L$ intersects ${\mathcal H}_r$ transversely and 
${\mathcal H}_r \cap L$ is isomorphic to ${\mathcal H}_{r-1}$ as
in the case of the hyperplane at infinity above,
or else $L$ is the tangent hyperplane to ${\mathcal H}_r$ at
an $\bF_{q^2}$-rational point $p$ and in that case
${\mathcal H}_r \cap L$ is isomorphic to the cone over
${\mathcal H}_{r-2}$ with vertex at $p$.  

If $r$ is even, using \eqref{nopts} above, it is easily 
checked that the largest number of affine
$\bF_{q^2}$-rational points on $L$ is attained when $L$ is
tangent to ${\mathcal H}_r$ at a point $p$ in $X_0 = 0$. 
This number is $z = q^{2r-1} + q^r$ and $d = n - z$ gives
the desired result.  On the other hand, if $r$ is odd, then
the maximum is attained when $L, {\mathcal H}_r,$ and the 
hyperplane $X_0$ intersect transversely.  We
have $z = q^{2r-1} + q^{r-1}$ in this case.
\end{proof}

For example, with $q = 2$, the Hermitian hypersurfaces yield 
$C_1$ codes as follows
over $\bF_4$:
\begin{align*}
          r = 2 &\qquad [n,k,d] = [36, 4, 24]\\
          r = 3 &\qquad [n,k,d] = [120, 5, 84]\\
          r = 4 &\qquad [n,k,d] = [528, 6, 384]
\end{align*}
For large $q$, in an asymptotic sense the $C_1$ codes come  
close to achieving the Griesmer bound on $n$ for the given
$k$ and $d$.  For instance, with $r = 5$, the Griesmer
bound gives 
$$n \ge q^{11} - q^5 - q^4 - q^3 - q^2 - q - 1$$
for a code with $k = r + 2 = 7$ and $d = q^{11} - q^9 - q^5 - q^4$
over $\bF_{q^2}$.   The actual $n$ for our $C_1$ code
is $q^{11} - q^5$.  The ratio of the bound and the actual
$n$ tends to $1$ as $q \to \infty$.  

The minimum distances of the $C_a$ codes for $a \ge 2$
can be estimated by the tools in \cite{11} and \cite{19}, but
determining the exact $d$ is somewhat subtle in these cases
because of the many cases that must be considered.

As noted in \cite{20} and \cite{11}, there is also 
a close connection between Hermitian hypersurface
codes and codes from the Deligne-Lusztig varieties
from one class of algebraic groups over finite fields.  
Deligne-Lusztig varieties are a class
of varieties with many rational points over finite fields. 
Indeed they often attain the maximum number for 
varieties with their invariants.
For instance, the Deligne-Lusztig variety of type 
${}^2A_2$ is just the Hermitian curve ${\mathcal H}_1$
in $\bP^2$, which has the maximum possible number
of points for a curve of its genus allowed by the 
Hasse-Weil bound.   The variety of type 
${}^2A_3$ is the blow-up of the Hermitian surface
${\mathcal H}_2$ at its $\bF_{q^2}$-rational points and 
is again maximal.
The variety of type ${}^2A_4$ is obtained from 
the complete intersection of the Hermitian
3-fold ${\mathcal H}_3$ in $\bP^4$ and the hypersurface
$$X_0^{q^3+1} + X_1^{q^3+1} + \cdots + X_4^{q^3+1} = 0,$$
a surface with singularities precisely at the $(q^5+1)(q^2+1)$
$\bF_{q^2}$-rational points of ${\mathcal H}_3$.  Each of those
singular points blows up to a Hermitian curve on the 
Deligne-Lusztig variety.   In each case, some of the more 
accessible Deligne-Lusztig codes ``come from'' $C_a$
evaluation codes on the Hermitian hypersuface (see 
Proposition 4.1 and Remark 4.1 of \cite{11}).  
 
To conclude this section, we note that there is a different 
way to construct order domains from the Hermitian hypersurfaces
which gives a direct parallel to another form of the 
equation of the Hermitian curve.  In the case $r = 1$, it is also 
common to consider a linear change
of coordinates (see \cite{21}, Example VI.4.3):
$$(X,Y,Z) = 
(\delta X_2,\delta(\gamma+1)X_0+\gamma X_1,\delta X_0 + X_1),$$
where $\delta,\gamma \in \bF_{q^2}$ satisfy 
$\delta^{q+1} = \gamma^q +\gamma= -1$.
In geometric terms, this has the effect of making the line
at infinity tangent to the curve ${\mathcal H}_1$
at an ${\bF}_{q^2}$-rational point.
The intersection multiplicity of the tangent line with the curve at the 
point of tangency is $q+1$.  The corresponding equation has the familiar form
$$X^{q+1} = Y^qZ + YZ^q.$$
{}From this, we can see easily that the hypotheses of Theorem~\ref{GP}
are satisfied if we define an order function on the affine coordinate
ring 
\begin{equation}
\label{stdherm}
R = \bF_{q^2}[X,Y]/\langle X^{q+1} - Y^q - Y\rangle
\end{equation}
by $\rho(X) = q$, $\rho(Y) = q+1$. The value semigroup is
$\Gamma = \langle q, q+1\rangle \subset \bZ_{\ge 0}$, and 
as in \S 4, we can view the usual order domain
associated to the Hermitian curve as a deformation
of the monomial algebra $\bF_{q^2}[t^q,t^{q+1}]$ (the coordinate ring
of a singular monomial curve).

\begin{example}
For the Hermitian surface, we can perform exactly the same 
change of coordinates used in the curve case in the variables 
$X_0,X_2,X_3$, leaving the other variable $X_1$ unchanged.  
(There is an analogous transformation, of course, for the 
Hermitian hypersurface of any dimension.)  

The result for the Hermitian surface, for instance,
is an equation of the form
$$X^{q+1} + X_1^{q+1} = Y^q Z + Y Z^q.$$
As compared with the first construction above, this form puts
a tangent plane to the surface rather than a plane meeting the 
surface transversely as the plane at infinity.  If we dehomogenize by
setting $Z = 1$,
then we claim that the corresponding affine algebra
\begin{equation}
\label{hsurf}
R' = \bF_{q^2}[X,X_1,Y]/\langle X^{q+1} + X_1^{q+1} - Y^q - Y\rangle
\end{equation}
has an order function $\rho$ defined by the matrix
$$M' = \begin{pmatrix} q & 0 & q+1\\
                       0 & 1 & 0  \\ \end{pmatrix}$$
(for example).  The $>_{M',lex}$ order makes the leading 
term $X^{q+1}$ and there are exactly two terms of maximal $M'$-weight.
The monomials in the complement of the initial ideal are
$$\Delta' = \{X^a X_1^b Y^c : 0 \le a \le q; b,c\ge 0\}.$$
An easy argument parallel to the one above shows that the 
monomials $\Delta'$ have distinct $M$-weights.  Hence
the hypotheses of Geil and Pellikaan's theorem are satisfied
here too, and we get an order function on $R'$.  
The corresponding valuation on the function field
of the Hermitian surface is not of the form used in 
Theorem~\ref{cdof}, however.  Indeed, the easiest way to describe
this order domain seems to be as a deformation of 
the product of order domains of the curve \eqref{stdherm} and a line
(see \cite{8}, Proposition 7.5).

By the facts used in the proof of Theorem~\ref{hermcodes}, 
the tangent planes to the Hermitian surface at all 
$\bF_{q^2}$-rational points intersect the surface ${\mathcal H}_2$
in reducible curves, each consisting of $q+1$ distinct lines
passing through the point of tangency.  This reduces the number
of affine $\bF_{q^2}$-rational points on the variety defined
by \eqref{hsurf} to $q^5$ (from $q^5 + q^2$
as in Theorem~\ref{hermcodes}).  But this loss may be compensated to an 
extent by some of the other properties of the resulting codes.

For example, we observe that all the evaluation and dual evaluation
codes constructed from $R'$ in \eqref{hsurf} have the same sort of 
quasicyclic structure studied in \cite{13}.  Namely, this affine 
piece of the Hermitian surface has a large automorphism group
containing a cyclic subgroup $H$ of order
$q^2-1$ generated by the automorphism
$$\sigma : (X,X_1,Y) \mapsto (\alpha X, \alpha X_1, \alpha^{q+1}Y)$$
where $\alpha$ is any primitive element of $\bF_{q^2}$.
$\diamondsuit$
\end{example}

\begin{proposition}
The $q^5$ affine $\bF_{q^2}$-rational points decompose into
$q^3+q$ orbits of size $q^2 - 1$, one orbit of size $q-1$
and one orbit of size $1$ under the action of $H$.
\end{proposition}

\begin{proof}
The point $(0,0,0)$ is clearly the only fixed point
of $H$.  The other $q-1$ points on the line $X = X_1 = 0$
form the orbit of size $q - 1$.  The orbit of each other
$\bF_{q^2}$-rational point has size $q^2-1$. 
\end{proof}

As in \cite{13}, this gives all of the evaluation and dual evaluation 
codes constructed from $R'$ the structure of modules over
$\bF_{q^2}[t]$, and compact representations for encoding.

\section{Order Domains from Grassmannians and Flag Varieties}

The codes from  Grassmannians and flag varieties that come
from the construction we will describe have been studied by 
Rodier in \cite{19}, but the connection with order domains
was not studied in that work.  Geil has constructed
order domains from the Grassmannians $G(2,n)$ (see below
for the notation) in \cite{7} (see Example 9.4 of \cite{8}).  
The new observation here is that 
{\it all} of these Grassmannians and flag varieties, and the
codes constructed from them, can be studied within the context of 
order domains.

Theorem~\ref{fldef} shows that any known explicit flat deformation
of a variety to a toric variety ({\it toric deformations} for 
short) can be used to produce examples of 
order domains, using the extrinsic characterization in Theorem~\ref{GP}.  
In this section, we will see that this is true in particular for 
Grassmannians, and for flag varieties.  The fact that toric 
deformations exist in these cases has been established in algebraic 
geometry (because of connections with combinatorics, mirror 
symmetry, and the theory of singularities).

The book \cite{22} is an excellent general reference 
for the aspects of the theory of toric varieties
of greatest relevance here.
We will use several results presented in Chapter 11 of \cite{22}
to construct our examples.  We will begin with the case of Grassmannians 
and furnish an answer to a question posed in \cite{8},
Example 9.4, concerning the varieties defined
by minors of arbitrary size of generic matrices. The 
Grassmannians correspond to the case of the maximal minors
of a rectangular matrix.  

Recall that the Grassmannian $G(k,n)$ is a projective variety 
whose points are in one-to-one correspondence with 
the $k$-dimensional vector subspaces of an
$n$-dimensional vector space (or the ($k-1$)-dimensional linear
subvarieties in ${\mathbb P}^{n-1}$).  We will not use the projective 
case here.  We very briefly recall the construction, since it 
shows that the $\bF_q$-rational points of $G(k,n)$
correspond to linear subspaces defined over $\bF_q$,
it explains the connection with minors of matrices, and it
also shows how to find toric deformations of Grassmannians.  

We write $\overline{\bF}$ for an algebraic closure of our field.
Given any basis $\{v_1,\ldots,v_k\}$ for a $k$-dimensional
vector subspace $W$ of $\overline{\bF}^n$, we form the 
$k\times n$ matrix with rows $v_i$.  The $k\times k$
(maximal) minors of this matrix are components of the 
{\it Pl\"ucker coordinate vector} of $W$ in 
${\mathbb P}^{\binom{n}{k} - 1}$.
This is a well-defined invariant of $W$
because a change of basis multiplies 
all components of the Pl\"ucker coordinate vector by 
a nonzero constant (the determinant of the change of 
basis matrix).  Hence any choice of basis in $W$
yields the same point in ${\mathbb P}^{\binom{n}{k} - 1}$.  
The locus of all such points (for all $W$) forms the 
Grassmannian $G(k,n)$.

Our extrinsic construction of order domains from 
Grassmannians relies on the following fact from \cite{22}.

\begin{proposition}[\cite{22}, Proposition 11.10]
\label{grass}
There exists a toric deformation taking $G(k,n)$ to the 
projective toric variety defined by the semigroup $\Gamma_{k,n}$
defined as follows:  Let $N = (t_{ij})$ be a generic
$k\times n$ matrix (the $t_{ij}$ are independent indeterminates).
Then $\Gamma_{k,n}$ is the semigroup in $\bZ_{\ge 0}^{\binom{n}{k}}$
generated by the columns of the 
$\binom{n}{k} \times kn$ matrix ${\mathcal B}_{k,n}$,
whose $\ell$th row has 1's in the positions corresponding
to the $t_{ij}$ in the {\it diagonal} of the $\ell$th minor of $N$, 
and zeroes in all other positions.  
\end{proposition}

This is deduced from properties of canonical (``SAGBI'') bases
for subalgebras of polynomial rings (like the monomial
algebra $\bF_q[\Gamma]$ above) in \cite{22}.  

\begin{corollary}
For all $k,n$, the Grassmannians $G(k,n)$
yield order domains with value semigroup $\Gamma_{k,n}$.
\end{corollary}

\begin{proof}
This follows immediately from Theorem~\ref{fldef} and Proposition~\ref{grass}.
\end{proof}

We will illustrate the conclusion by considering the case of 
$G(3,5)$ here.  This is one of the simplest cases not covered by 
results in \cite{7} quoted in Example 9.4 of \cite{8}, and shows how the 
general proof of 
Theorem~\ref{grass} works.

\begin{example}
Let 
$$N = \begin{pmatrix} t_{11} &t_{12} &\cdots &t_{15}\\
               t_{21} &t_{22} &\cdots &t_{25}\\
               t_{31} &t_{32} &\cdots &t_{35}\\
\end{pmatrix}$$
be the generic $3\times 5$ matrix.  There are $\binom{5}{3} = 10$
maximal minors of $N$.  The diagonal terms are
\begin{equation}
\label{diags}
t_{11}t_{22}t_{33}, t_{11}t_{22}t_{34}, \ldots, t_{13}t_{24}t_{35}.
\end{equation}
The corresponding matrix ${\mathcal B}_{3,5}$ as in the statement of 
Proposition~\ref{grass} is
the following (columns are indexed by the entries of $N$, listed
row-wise):
\setcounter{MaxMatrixCols}{15}
$${\mathcal B}_{3,5}=
\begin{pmatrix} 
1 & 0 & 0 & 0 & 0 & 0 & 1 & 0 & 0 & 0 & 0 & 0 & 1 & 0 & 0\\
1 & 0 & 0 & 0 & 0 & 0 & 1 & 0 & 0 & 0 & 0 & 0 & 0 & 1 & 0\\
1 & 0 & 0 & 0 & 0 & 0 & 1 & 0 & 0 & 0 & 0 & 0 & 0 & 0 & 1\\
1 & 0 & 0 & 0 & 0 & 0 & 0 & 1 & 0 & 0 & 0 & 0 & 0 & 1 & 0\\
1 & 0 & 0 & 0 & 0 & 0 & 0 & 1 & 0 & 0 & 0 & 0 & 0 & 0 & 1\\
1 & 0 & 0 & 0 & 0 & 0 & 0 & 0 & 1 & 0 & 0 & 0 & 0 & 0 & 1\\
0 & 1 & 0 & 0 & 0 & 0 & 0 & 1 & 0 & 0 & 0 & 0 & 0 & 1 & 0\\
0 & 1 & 0 & 0 & 0 & 0 & 0 & 1 & 0 & 0 & 0 & 0 & 0 & 0 & 1\\
0 & 1 & 0 & 0 & 0 & 0 & 0 & 0 & 1 & 0 & 0 & 0 & 0 & 0 & 1\\
0 & 0 & 1 & 0 & 0 & 0 & 0 & 0 & 1 & 0 & 0 & 0 & 0 & 0 & 1\\
\end{pmatrix} .
$$

The columns of ${\mathcal B}_{3,5}$ generate the semigroup $\Gamma_{3,5}$.
Write $X_1,\ldots,X_{10}$ for the coordinates in ${\mathbb P}^{10}$.
The toric variety corresponding to $\bF_q[\Gamma_{3,5}]$ is given by the 
parametrization
$$X_1 = t_{11}t_{22}t_{33}, X_2 = t_{11}t_{22}t_{34}, \ldots, 
X_{10} = t_{13}t_{24}t_{35}$$
(using \eqref{diags}).  Eliminating the $t_{ij}$, we find the graded 
reverse lex Gr\"obner basis of $I_{\Gamma_{3,5}}$ equals
\begin{align}
\label{tgbasis}
G_T = \{ X_8 X_6 &- X_9 X_5,\notag\\
         X_7 X_6 &- X_4 X_9,\notag\\ 
         X_7 X_5 &- X_8 X_4,\\
         X_7 X_3 &- X_8 X_2,\notag\\
         X_4 X_3 &- X_5 X_2\}\notag
\end{align}
(the positive term is the leading term in each case).
The corresponding projective toric variety has dimension 6 and degree 5 in 
${\mathbb P}^9$; the equations above can also be viewed as defining
the affine cone over that projective variety, which has 
dimension 7 in ${\mathbb A}^{10}$.

The ideal of the Grassmannian $G(3,5)$ is generated by quadratic
polynomials called the {\it Pl\"ucker relations} between the 
Pl\"ucker coordinate vectors of 3-planes $W$.  We have the following
Gr\"obner basis for this ideal with respect to the same graded
reverse lex order as in \eqref{tgbasis}:
\begin{align}
\label{ggbasis}
G_G = 
\{X_8 X_6 &- X_9 X_5 + X_3 X_{10} , \notag\\
    X_7 X_6 &- X_4 X_9 + X_2 X_{10}, \notag\\
    X_7 X_5 &- X_8 X_4 + X_1 X_{10}, \\
    X_7 X_3 &- X_8 X_2 + X_1 X_9,\notag\\
    X_4 X_3 &- X_5 X_2 + X_1 X_6.\}\notag
\end{align}

To understand the order domain structure here, we need to introduce
the weight matrix $M = {\mathcal B}_{3,5}^t$, a $15\times 10$ matrix
of rank $7$.  In our discussion of the weight orders in Geil
and Pellikaan's theorem (Theorem~\ref{GP}), note that we used only 
matrices where the number of rows equals the rank of the corresponding 
valuation
of the function field.  That is not necessary, though.  It would 
be perfectly legal to define a matrix weight order using the 
full matrix $M$;  the value semigroup is then a sub-semigroup 
in $\bZ_{\ge 0}^{15}$ whose rank is 7.

Note that in each polynomial in \eqref{ggbasis}, the same two terms 
as in the corresponding polynomial in \eqref{tgbasis} appear.  These 
are the terms of maximum $M$-weight in each case.  The remaining
terms in \eqref{ggbasis} have smaller $M$-weight.  So the ideal in 
\eqref{tgbasis} is indeed a toric deformation of the ideal in 
\eqref{ggbasis}.

Moreover by the construction here, two monomials in the $X_i$, say
$X^\alpha$ and $X^\beta$, have the same weight if and only 
if $M\alpha = M\beta$.  But that implies $X^\alpha - X^\beta$
is in the toric ideal $I_{\Gamma_{3,5}}$, so one of the monomials
is divisible by one of the leading terms of $G_G$ or $G_T$.
This shows that \eqref{ggbasis} satisfies all the hypotheses of 
Theorem~\ref{GP}.  Hence we have constructed an order domain from the 
Grassmannian $G(3,5)$ (or more properly, the affine cone over the 
Grassmannian) with the homogeneous ideal given by \eqref{ggbasis}.
We could also set $X_1 = 1$ to obtain an affine algebra
corresponding to an affine subset of the 
Grassmannian itself if we wish, and the corresponding weight 
matrix $M'$ would be the rank 6 matrix obtained by transposing 
the submatrix of ${\mathcal B}_{3,5}$ formed by omitting the first 
row.  
$\diamondsuit$
\end{example}

Just as the Grassmannian $G(k,n)$ is a 
projective variety whose points correspond to the $k$-dimensional vector
subspaces of $n$-space, the partial flag variety $F(n_1,\ldots,n_\ell;n)$
is a variety whose points correspond to partial flags of vector subspaces:
$$V_1 \subset V_2 \subset \cdots \subset V_{\ell}$$
in $n$-space, where $\dim(V_i) = n_i$ for all $i = 1, \ldots, \ell$.
By considering the Pl\"ucker embeddings, we have a natural inclusion
$$F(n_1,\ldots,n_\ell;n) \subset {\mathbb P}^{N_1-1} 
\times \cdots \times {\mathbb P}^{N_\ell - 1},$$
where $N_i = \binom{n}{n_i}$.  The flag variety is then defined
in this product of projective spaces by the conditions that
$V_1 \subset V_2 \subset \cdots V_\ell$.  If desired, the 
product of projective spaces can also be embedded in a single
projective space by the usual Segre mapping.

The existence of toric deformations of all $F(n_1,\ldots,n_\ell;n)$
is known from work of Gonciulea and Lakshmibai (\cite{9}, Theorem 10.6).  
A more combinatorial description similar to that provided by 
Sturmfels for the Grassmannians has appeared in \cite{15}.
We review this and indicate how to derive explicit order domains from 
$F(n_1,\ldots,n_\ell;n)$.  

The flag variety $F(n_1,\ldots,n_\ell;n)$ can be identified
with the quotient $SL(n)/Q$, where $Q = \cap_{i=1}^\ell P_{n_i}$,
and $P_{n_i}$ is the parabolic subgroup of $SL(n)$ consisting
of matrices of the form 
$\bigl( \begin{smallmatrix} * & *\\ 0 & *\\\end{smallmatrix}\bigr)$
with $0$ an $(n-n_i)\times n_i$ zero matrix.  Let 
$H = \cup_{i=1}^\ell W^{n_i}$, where
$$W^{n_i} = \{(j_1,\ldots,j_{n_i}): 1\le j_1 < j_2 < \cdots 
< j_{n_i}\le n\}.$$
Note that the elements of $H$ are in one-to-one correspondence with the 
Pl\"ucker coordinates on 
$$G(n_1,n) \times G(n_2,n) \times \cdots \times G(n_\ell,n).$$ 
Given $\pi \in H$, we will write $p_\pi$ for the corresponding
Pl\"ucker coordinate. 

The following construction defines a partial order $\succ$ on $H$.  Let
$\pi = (i_1,\ldots,i_a)$ and $\pi'=(j_1,\ldots,j_b)$
in $H$.  Let
$$\pi \succeq \pi' \Leftrightarrow a \le b \hbox{ and } i_s \ge j_s, 
s=1,\ldots,a.$$
The set $H$ is a {\it distributive lattice} under the partial order 
relation $\succeq$. 

For each pair of elements $\pi,\pi'$ that are incomparable in 
the ordering on $H$, there
is a quadratic (Grassmann-Pl\"ucker) relation
\begin{equation}
\label{GrassPluck}
p_\pi p_{\pi'} = \sum c_{\lambda \mu} p_\lambda p_\mu.
\end{equation}
Gonciulea and Lakshmibai show that the monomial 
$p_{\max(\pi,\pi')}p_{\min(\pi,\pi')}$ appears
with coefficient 1 on the right hand side of \eqref{GrassPluck} and 
the other terms are smaller with respect to a suitable
monomial order.  Moreover, the relations \eqref{GrassPluck} are 
a Gr\"obner basis for the ideal of $F(n_1,\ldots,n_\ell;n)$
in ${\mathbb P}^{N_1-1}\times \cdots\times {\mathbb P}^{N_\ell - 1}$.
The equations
$$p_\pi p_{\pi'} - p_{\max(\pi,\pi')}p_{\min(\pi,\pi')} = 0$$
define a toric subvariety $X = X(n_1,\ldots,n_\ell;n)$
of $G(n_1,n)\times \cdots G(n_\ell,n)$
and the flag variety has a flat deformation to $X$.

\begin{example}
We consider the flag varieties $F(1,n-1;n)$ studied 
in \cite{19}.  In particular, that 
article shows the corresponding codes compare very favorably
with projective Reed-Muller codes.  The same techniques would be applicable 
to all these flag varieties.  In this case the set $H$ defined
above reduces to 
$$H = \{(1),(2),\ldots,(n),(\hat{1}),(\hat{2}),\ldots,(\hat{n})\},$$
where $\hat{j} = (1,\ldots,j-1,j+1,\ldots,n)$.  Using the definition
of the partial order $\succeq$, we have
$$(n) \succeq (n-1) \succeq \cdots \succeq (1),$$
and 
$$(2) \succeq (\hat{1}) \succeq \cdots \succeq (\hat{n}).$$
There is exactly one pair of incomparable elements:
$(1)$ and $(\hat{1})$, and $\max((1),(\hat{1})) = (2)$,
$\min((1),(\hat{1})) = (\hat{2})$.  The corresponding relation 
\eqref{GrassPluck} is
$$p_{(1)} p_{(\hat{1})} = p_{(2)}p_{(\hat{2})} - p_{(3)} p_{(\hat{3})}
+ \cdots + (-1)^{n-1} p_{(n)}p_{(\hat{n})}.$$
(This equation, obtained from the obvious determinant 
expansion, expresses the condition $V_1 \subset V_2$
where $\dim(V_1) = 1$ and $\dim(V_2) = n-1$.)
The toric deformation is defined by
$$p_{(1)} p_{(\hat{1})} = p_{(2)}p_{(\hat{2})}$$
in the product ${\mathbb P}^{n-1}\times {\mathbb P}^{n-1}$.
We could also embed $F(1,n-1;n)$ in ${\mathbb P}^{n^2-1}$ via
the standard Segre mapping.  There seems to be 
little advantage in doing that, however, because of the 
large number of additional equations needed to define the 
ideal of the Segre image.
$\diamondsuit$
\end{example}




\end{document}